\documentclass[12pt,twoside]{article}
\usepackage{amsmath,amssymb}
\usepackage{amsthm}
\usepackage{upref}
\usepackage{citeref}
\numberwithin{equation}{section}

\theoremstyle{plain}
\newtheorem{theorem}{Theorem}[section]
\newtheorem{proposition}[theorem]{Proposition}
\newtheorem{definition}[theorem]{Definition}
\newtheorem{corollary}[theorem]{Corollary}

\theoremstyle{definition}
\newtheorem{remark}[theorem]{Remark}

\newtheorem{example}[theorem]{Example}

\begin{document}

\newcommand{\eq}{equation}
\newcommand{\real}{\ensuremath{\mathbb R}}
\newcommand{\comp}{\ensuremath{\mathbb C}}
\newcommand{\rn}{\ensuremath{{\mathbb R}^n}}
\newcommand{\tn}{\ensuremath{{\mathbb T}^n}}
\newcommand{\rnp}{\ensuremath{\real^n_+}}
\newcommand{\rnn}{\ensuremath{\real^n_-}}
\newcommand{\Rn}{\ensuremath{{\mathbb R}^{n-1}}}
\newcommand{\no}{\ensuremath{\nat_0}}
\newcommand{\ganz}{\ensuremath{\mathbb Z}}
\newcommand{\zn}{\ensuremath{{\mathbb Z}^n}}
\newcommand{\zom}{\ensuremath{{\mathbb Z}_{\Om}}}
\newcommand{\zOm}{\ensuremath{{\mathbb Z}^{\Om}}}
\newcommand{\As}{\ensuremath{A^s_{p,q}}}
\newcommand{\Bs}{\ensuremath{B^s_{p,q}}}
\newcommand{\Fs}{\ensuremath{F^s_{p,q}}}
\newcommand{\Fsr}{\ensuremath{F^{s,\rloc}_{p,q}}}
\newcommand{\nat}{\ensuremath{\mathbb N}}
\newcommand{\Om}{\ensuremath{\Omega}}
\newcommand{\di}{\ensuremath{{\mathrm d}}}
\newcommand{\sn}{\ensuremath{{\mathbb S}^{n-1}}}
\newcommand{\Ac}{\ensuremath{\mathcal A}}
\newcommand{\Acs}{\ensuremath{\Ac^s_{p,q}}}
\newcommand{\Bc}{\ensuremath{\mathcal B}}
\newcommand{\Cc}{\ensuremath{\mathcal C}}
\newcommand{\cc}{{\scriptsize $\Cc$}${}^s (\rn)$}
\newcommand{\ccd}{{\scriptsize $\Cc$}${}^s (\rn, \delta)$}
\newcommand{\Fc}{\ensuremath{\mathcal F}}
\newcommand{\Lc}{\ensuremath{\mathcal L}}
\newcommand{\Mc}{\ensuremath{\mathcal M}}
\newcommand{\Ec}{\ensuremath{\mathcal E}}
\newcommand{\Pc}{\ensuremath{\mathcal P}}
\newcommand{\Efr}{\ensuremath{\mathfrak E}}
\newcommand{\Mfr}{\ensuremath{\mathfrak M}}
\newcommand{\Mbf}{\ensuremath{\mathbf M}}
\newcommand{\Dbb}{\ensuremath{\mathbb D}}
\newcommand{\Lbb}{\ensuremath{\mathbb L}}
\newcommand{\Pbb}{\ensuremath{\mathbb P}}
\newcommand{\Qbb}{\ensuremath{\mathbb Q}}
\newcommand{\Rbb}{\ensuremath{\mathbb R}}
\newcommand{\vp}{\ensuremath{\varphi}}
\newcommand{\hra}{\ensuremath{\hookrightarrow}}
\newcommand{\supp}{\ensuremath{\mathrm{supp \,}}}
\newcommand{\ssupp}{\ensuremath{\mathrm{sing \ supp\,}}}
\newcommand{\dist}{\ensuremath{\mathrm{dist \,}}}
\newcommand{\unif}{\ensuremath{\mathrm{unif}}}
\newcommand{\ve}{\ensuremath{\varepsilon}}
\newcommand{\vk}{\ensuremath{\varkappa}}
\newcommand{\vr}{\ensuremath{\varrho}}
\newcommand{\pa}{\ensuremath{\partial}}
\newcommand{\oa}{\ensuremath{\overline{a}}}
\newcommand{\ob}{\ensuremath{\overline{b}}}
\newcommand{\of}{\ensuremath{\overline{f}}}
\newcommand{\LA}{\ensuremath{L^r\!\As}}
\newcommand{\LcA}{\ensuremath{\Lc^{r}\!A^s_{p,q}}}
\newcommand{\LcdA}{\ensuremath{\Lc^{r}\!A^{s+d}_{p,q}}}
\newcommand{\LcB}{\ensuremath{\Lc^{r}\!B^s_{p,q}}}
\newcommand{\LcF}{\ensuremath{\Lc^{r}\!F^s_{p,q}}}
\newcommand{\Lb}{\ensuremath{L^r b^s_{p,q}}}
\newcommand{\Lf}{\ensuremath{L^r\!f^s_{p,q}}}
\newcommand{\La}{\ensuremath{L^r a^s_{p,q}}}
\newcommand{\Lob}{\ensuremath{L^r \ob{}^s_{p,q}}}
\newcommand{\Lof}{\ensuremath{L^r \of{}^s_{p,q}}}
\newcommand{\Loa}{\ensuremath{L^r\, \oa{}^s_{p,q}}}
\newcommand{\Lcoa}{\ensuremath{\Lc^{r}\oa{}^s_{p,q}}}
\newcommand{\Lcob}{\ensuremath{\Lc^{r}\ob{}^s_{p,q}}}
\newcommand{\Lcof}{\ensuremath{\Lc^{r}\of{}^s_{p,q}}}
\newcommand{\Lca}{\ensuremath{\Lc^{r}\!a^s_{p,q}}}
\newcommand{\Lcb}{\ensuremath{\Lc^{r}\!b^s_{p,q}}}
\newcommand{\Lcf}{\ensuremath{\Lc^{r}\!f^s_{p,q}}}
\newcommand{\id}{\ensuremath{\mathrm{id}}}
\newcommand{\tr}{\ensuremath{\mathrm{tr\,}}}
\newcommand{\trd}{\ensuremath{\mathrm{tr}_d}}
\newcommand{\trL}{\ensuremath{\mathrm{tr}_L}}
\newcommand{\ext}{\ensuremath{\mathrm{ext}}}
\newcommand{\re}{\ensuremath{\mathrm{re\,}}}
\newcommand{\Rea}{\ensuremath{\mathrm{Re\,}}}
\newcommand{\Ima}{\ensuremath{\mathrm{Im\,}}}
\newcommand{\loc}{\ensuremath{\mathrm{loc}}}
\newcommand{\rloc}{\ensuremath{\mathrm{rloc}}}
\newcommand{\osc}{\ensuremath{\mathrm{osc}}}
\newcommand{\pr}{\pageref}
\newcommand{\wh}{\ensuremath{\widehat}}
\newcommand{\wt}{\ensuremath{\widetilde}}
\newcommand{\ol}{\ensuremath{\overline}}
\newcommand{\os}{\ensuremath{\overset}}
\newcommand{\Li}{\ensuremath{\overset{\circ}{L}}}
\newcommand{\Ai}{\ensuremath{\os{\, \circ}{A}}}
\newcommand{\Ci}{\ensuremath{\os{\circ}{\Cc}}}
\newcommand{\dom}{\ensuremath{\mathrm{dom \,}}}
\newcommand{\SA}{\ensuremath{S^r_{p,q} A}}
\newcommand{\SB}{\ensuremath{S^r_{p,q} B}}
\newcommand{\SF}{\ensuremath{S^r_{p,q} F}}
\newcommand{\Hc}{\ensuremath{\mathcal H}}
\newcommand{\Lci}{\ensuremath{\overset{\circ}{\Lc}}}
\newcommand{\bmo}{\ensuremath{\mathrm{bmo}}}
\newcommand{\BMO}{\ensuremath{\mathrm{BMO}}}
\newcommand{\cm}{\\[0.1cm]}
\newcommand{\Aa}{\ensuremath{\os{\, \ast}{A}}}
\newcommand{\Ba}{\ensuremath{\os{\, \ast}{B}}}
\newcommand{\Fa}{\ensuremath{\os{\, \ast}{F}}}
\newcommand{\Aas}{\ensuremath{\Aa{}^s_{p,q}}}
\newcommand{\Bas}{\ensuremath{\Ba{}^s_{p,q}}}
\newcommand{\Fas}{\ensuremath{\Fa{}^s_{p,q}}}
\newcommand{\Ca}{\ensuremath{\os{\, \ast}{{\mathcal C}}}}
\newcommand{\Cas}{\ensuremath{\Ca{}^s}}
\newcommand{\Car}{\ensuremath{\Ca{}^r}}

\begin{center}
{\Large Truncations and compositions in function spaces}
\\[1cm]
{Hans Triebel}
\\[0.2cm]
Institut f\"{u}r Mathematik\\
Friedrich--Schiller--Universit\"{a}t Jena\\
07737 Jena, Germany
\\[0.1cm]
email: hans.triebel@uni-jena.de
\\[5cm]
\end{center}

\begin{abstract}
These notes deal with some recent assertions about truncations $f \mapsto |f|$ and compositions $f \mapsto g\circ f$ in the spaces 
$\As (\rn)$, $A \in \{B,F \}$.
\\[2cm]
\end{abstract}



\newpage

\tableofcontents

\newpage

\section{Introduction}   \label{S1}
Let $\As (\rn)$, $A\in \{B,F \}$ with $0<p,q \le \infty$ and $s\in \real$ be the nowadays well--known function spaces in Euclidean 
$n$--space. If
\begin{\eq}   \label{1.1}
s> \sigma^n_p = n \Big( \max \big(\frac{1}{p},1 \big) - 1 \Big), \qquad n\in \nat, \quad \text{and} \quad 0<p \le \infty,
\end{\eq}
then $\As (\rn)$ consists entirely  of regular tempered distributions, $\As (\rn) \subset L^{\loc}_1 (\rn)$. Under these circumstances 
it makes sense to ask whether $|f| \in \As (\rn)$ if $f\in \As (\rn)$ real and
\begin{\eq}   \label{1.2}
\big\| \, |f| \, |\As (\rn) \big\| \le  c \, \|f \, | \As (\rn) \|, \qquad f \in \As (\rn) \quad \text{real},
\end{\eq}
for some $c>0$ and all $f\in \As (\rn)$ real. This is the {\em truncation property} which attracted  some attention since decades. We 
say that $\As (\rn)$ has the {\em strong truncation property} if \eqref{1.2} can be strengthened  by the equivalence
\begin{\eq}   \label{1.3}
\big\| \, |f| \, |\As (\rn) \big\| \sim  \, \|f \, | \As (\rn) \|, \qquad f \in \As (\rn) \quad \text{real}.
\end{\eq}
The {\em perfect truncation property} asks for conditions ensuring that any real bounded continuous function $f$ belongs to $\As (\rn)$
if, and only if, $|f| \in \As (\rn)$, including \eqref{1.3}. Whereas truncations of type \eqref{1.2} have some history, the indicated
reinforcements in terms of strong and perfect truncations came up only very recently, \cite{MiVS19, MuN19, VScha20}. We deal with 
these problems in Section \ref{S3}. Let $g$ be a real and continuous function on $\real$ with $g(0) =0$. Let $\As (\rn)$ be as above
with $s > \sigma^n_p$ as in \eqref{1.1}. Then it makes sense to ask whether
\begin{\eq}    \label{1.4}
g \circ f = g \big( f(\cdot) \big) \in \As (\rn) \qquad \text{for any real $f\in \As (\rn)$}.
\end{\eq}
This so--called {\em composition} (of $f$ with $g$) attracted a lot of attention since decades. The most distinguished case of 
composition is the above truncation based on $g(t) =|t|$ for $t \in \real$. Furthermore  it is of interest to find conditions ensuring
that there is an $g$--version  of the {\em sublinearity} \eqref{1.2} which means that
\begin{\eq}   \label{1.5}
\| g \circ f \, | \As (\rn) \| \le c_g \, \|f \, | \As (\rn) \|, \qquad f \in \As (\rn) \quad \text{real}.
\end{\eq}
Afterwards it makes sense  to ask for counterparts of the above strong truncation and perfect truncation, resulting in {\em
strong $g$--composition} and {\em perfect $g$--composition}. This is the topic of Section \ref{S4}. In Section \ref{S2}\ we collect
definitions and a few distinguished ingredients for the spaces $\As (\rn)$, adapted to their use in connection with truncations and
compositions. Some of them might be of self--contained  interest, complementing of what is already known. We will indicate how these
specified  properties can be used in the context of truncations and compositions. But otherwise  one must consult the related sections 
in the collection \cite{Tri21} for proofs. In other words, we survey some rather specific recent aspects of truncations and 
compositions complementing what is already known.

\section{Function spaces}    \label{S2}
\subsection{Definitions}    \label{S2.1}
We use standard notation. Let $\nat$ be the collection of all natural numbers and $\no = \nat \cup \{0 \}$. Let $\rn$ be Euclidean $n$-space, where $n \in \nat$. Put $\real = \real^1$, whereas $\comp$ is the complex plane. Let $S(\rn)$ be the Schwartz space of all complex-valued rapidly decreasing
infinitely differentiable functions on $\rn$ and let $S' (\rn)$ be the  dual space of all
tempered distributions on \rn. Let $L_p (\rn)$ with $0<p \le \infty$ be the standard complex quasi-Banach space with respect to the Lebesgue measure in \rn, quasi-normed by
\begin{\eq}   \label{2.1}
\| f \, | L_p (\rn) \| = \Big( \int_{\rn} |f(x)|^p \, \di x \Big)^{1/p}
\end{\eq}
with the natural modification if $p= \infty$. Similarly $L_p (M)$ where $M$ is a Lebesgue--measurable subset of \rn.
As usual $\ganz$ is the collection of all integers; and $\zn$ where $n\in \nat$ denotes the lattice of all points $m = (m_1, \ldots, m_n) \in \rn$ with $m_j \in \ganz$. 
If $\vp \in S(\rn)$ then
\begin{\eq}   \label{2.2}
\wh{\vp} (\xi) = (F \vp)(\xi) = (2 \pi)^{-n/2} \int_{\rn} e^{-ix \xi} \, \vp (x) \, \di  x, \qquad \xi \in \rn,
\end{\eq}
denotes the Fourier transform of \vp. As usual, $F^{-1} \vp$ and $\vp^\vee$ stand for the inverse Fourier transform, which is given by the right-hand side of 
\eqref{2.2} with $i$ in place $-i$. Note that $x \xi = \sum^n_{j=1} x_j \xi_j$ 
is  the scalar product in \rn. Both $F$ and $F^{-1}$ are extended to $S' (\rn)$ in the standard way. Let $\vp_0 \in S(\rn)$ with
\begin{\eq}   \label{2.3}
\vp_0 (x) =1 \ \text{if $|x| \le 1$} \qquad \text{and} \qquad \vp_0 (x) =0 \ \text{if $|x| \ge 3/2$},
\end{\eq}
and let
\begin{\eq}   \label{2.4}
\vp_k (x) = \vp_0 \big( 2^{-k} x) - \vp_0 \big( 2^{-k+1} x \big), \qquad x \in \rn, \quad k \in \nat.
\end{\eq}
Since
\begin{\eq}    \label{2.5}
\sum^\infty_{j=0} \vp_j (x) = 1 \qquad \text{for} \quad x\in \rn,
\end{\eq}
the $\vp_j$ form a dyadic resolution of unity. The entire analytic functions $(\vp_j \wh{f} )^\vee (x)$ make sense pointwise in $\rn$ for any $f \in S' (\rn)$. Let
\begin{\eq}   \label{2.6}
Q_{J,M} = 2^{-J} M + 2^{-J} (0,1)^n, \qquad J \in \ganz, \quad M \in \zn,
\end{\eq}
be the usual dyadic cubes in \rn, $n\in \nat$, with sides of length $2^{-J}$ parallel to the axes of coordinates and with $2^{-J}M$
as lower left corner. If $Q$ is a cube in $\rn$ and $d>0$ then $dQ$ is the cube in $\rn$ concentric with $Q$ whose  side--length
is $d$ times the side--length of $Q$. Let $|\Om|$ be the Lebesgue measure of the Lebesgue measurable  set $\Om$ in \rn. Let $a^+ =
\max (a,0)$ for $a\in \real$.

\begin{definition}   \label{D2.1}
Let $\vp = \{ \vp_j \}^\infty_{j=0}$ be the above dyadic resolution  of unity.
\\[0.1cm]
{\upshape (i)} Let
\begin{\eq}   \label{2.7}
0<p \le \infty, \qquad 0<q \le \infty, \qquad s \in \real.
\end{\eq}
Then $\Bs (\rn)$ is the collection of all $f \in S' (\rn)$ such that
\begin{\eq}   \label{2.8}
\| f \, | \Bs (\rn) \|_{\vp} = \Big( \sum^\infty_{j=0} 2^{jsq} \big\| (\vp_j \widehat{f})^\vee \, | L_p (\rn) \big\|^q \Big)^{1/q} 
\end{\eq}
is finite $($with the usual modification if $q= \infty)$. 
\\[0.1cm]
{\upshape (ii)} Let
\begin{\eq}   \label{2.9}
0<p<\infty, \qquad 0<q \le \infty, \qquad s \in \real.
\end{\eq}
Then $F^s_{p,q} (\rn)$ is the collection of all $f \in S' (\rn)$ such that
\begin{\eq}   \label{2.10}
\| f \, | F^s_{p,q} (\rn) \|_{\vp} = \Big\| \Big( \sum^\infty_{j=0} 2^{jsq} \big| (\vp_j \wh{f})^\vee (\cdot) \big|^q \Big)^{1/q} \big| L_p (\rn) \Big\|
\end{\eq}
is finite $($with the usual modification if $q=\infty)$.
\\[0.1cm]
{\upshape (iii)} Let $0<q < \infty$ and $s\in \real$. Then $F^s_{\infty,q} (\rn)$ is the collection of all $f\in S'(\rn)$ such that
\begin{\eq}  \label{2.11}
\| f \, | F^s_{\infty,q} (\rn) \|_{\vp} = \sup_{J\in \ganz, M\in \zn} 2^{Jn/q} \Big(\int_{Q_{J,M}} \sum_{j \ge J^+} 2^{jsq} \big|
(\vp_j \wh{f} )^\vee (x) \big|^q \, \di x \Big)^{1/q}
\end{\eq}
is finite. Let $F^s_{\infty, \infty} (\rn) = B^s_{\infty, \infty} (\rn)$.
\end{definition}

\begin{remark}    \label{R2.2}
We take for granted that the reader is familiar with the theory of the spaces $\As (\rn)$, $A\in \{B,F \}$, $s\in \real$ and $0 <p,q
\le \infty$. In particular they are independent of the chosen resolution of unity $\vp$ (equivalent quasi--norms). This justifies our
omission of the subscript $\vp$ in \eqref{2.8}, \eqref{2.10} and \eqref{2.11}. Discussions and references may be found in 
\cite[Section 1.1.1, pp.\,1--5]{T20}. The spaces $F^s_{\infty, q} (\rn)$  will play only a marginal role. On the other hand it will be 
occasionally useful to specify some assertions to the following  special spaces. Let
\begin{\eq}   \label{2.12}
B^s_p (\rn) = B^s_{p,p} (\rn) = F^s_{p,p} (\rn), \quad 0<p \le \infty, \quad \text{and} \quad s\in \real,
\end{\eq}
including the {\em H\"{o}lder--Zygmund spaces}
\begin{\eq}   \label{2.13}
\Cc^s (\rn) = B^s_\infty (\rn) = B^s_{\infty,\infty} (\rn), \qquad s \in \real.
\end{\eq}
Furthermore,
\begin{\eq}   \label{2.14}
H^s_p (\rn) = F^s_{p,2} (\rn), \qquad 1<p<\infty \quad \text{and} \quad s\in \real
\end{\eq}
are the {\em Sobolev spaces}, including the {\em classical Sobolev spaces}
\begin{\eq}   \label{2.15}
W^k_p (\rn) = H^k_p (\rn), \qquad 1<p<\infty, \quad k\in \no,
\end{\eq}
which can be equivalently normed by
\begin{\eq}   \label{2.16}
\| f \, | W^k_p(\rn) \| = \sum_{|\alpha| \le k} \| D^\alpha f \, | L_p (\rn) \|, \qquad 1<p<\infty, \quad k\in \no.
\end{\eq}
\end{remark}

\subsection{Fubini property}   \label{S2.2}
In what follows we collect some specific properties of the spaces $\As (\rn)$ which will be useful in connection with truncations and
compositions as described in the Introduction so far. Some of them will be taken over from the literature, others will be adapted and
complemented.

Let $n\in \nat$ and
\begin{\eq}   \label{2.17}
\sigma^n_p = n \Big( \max \big( \frac{1}{p}, 1 \big) - 1 \Big) \quad \text{and} \quad \sigma^n_{p,q} = n \Big( \max \big(\frac{1}{p}, \frac{1}{q}, 1 \big) - 1 \Big),
\end{\eq}
where $0<p,q \le \infty$. Recall that $L^{\loc}_1 (\rn)$ consists of all Lebesgue--measurable  (complex--valued) functions in $\rn$
which are Lebesgue--integrable on any bounded domains (= open set) in \rn. Based on related references we characterized in \cite[Section
2.1.2, pp.\,23--24]{T20} which spaces consist entirely of regular distributions,
\begin{\eq}  \label{2.18}
\As (\rn) \subset L^{\loc}_1 (\rn).
\end{\eq}
This applies to all spaces
\begin{\eq}   \label{2.19}
\As (\rn) \qquad \text{with $0<p,q \le \infty$ and $s> \sigma^n_p$},
\end{\eq}
and some limiting spaces with $s= \sigma^n_p$. Then it makes sense what follows. Let $n \ge 2$,
\begin{\eq}   \label{2.20}
x = (x_1, \ldots, x_n ) \in \rn \quad \text{and} \quad x^j = (x_1,..,x_{j-1}, x_{j+1},..,x_n) \in \Rn.
\end{\eq}
Let $f \in L^{\loc}_1 (\rn)$. Then
\begin{\eq}   \label{2.21}
x_j \to f^{x^j} (x_j) = f(x), \qquad x\in \rn,
\end{\eq}
is considered as a function on $\real$ for any $x^j \in \Rn$ a.e. 

\begin{definition}    \label{D2.3}
Let $n \ge 2$. Let $0<p,q\le \infty$ and $s > \sigma^n_p$ with $p<\infty$ for the $F$--spaces. Then $\As (\rn)$ is said to have the
Fubini  property if for any $f\in \As (\rn)$ the function $f^{x^j}$ according to \eqref{2.21}  belongs to $\As (\real)$ for almost all
$x^j \in \Rn$ and
\begin{\eq}   \label{2.22}
\| f \, | \As (\rn) \| \sim \sum^n_{j=1} \big\| \, \| f^{x^j} \, | \As(\real) \| \, \big| L_p (\Rn) \big\|
\end{\eq}
$($equivalent quasi--norms$)$.
\end{definition}

\begin{remark}   \label{R2.4}
This coincides with \cite[Definition 4.2, p.\,35]{T01}. The modified version in \cite[Definition3.23, p.\,104]{T20} covers also the
spaces $F^s_{\infty, q} (\rn)$.
\end{remark}

\begin{proposition}   \label{P2.5}
Let $n \ge 2$.
\cm
{\em (i)} The spaces $\Fs (\rn)$ with
\begin{\eq}    \label{2.23}
0<p<\infty, \quad 0<q \le \infty \quad \text{and} \quad s > \sigma^n_{p,q},
\end{\eq}
have the Fubini property.
\cm
{\em (ii)} The spaces $\Bs (\rn)$ with $0<p,q \le \infty$ and $s> \sigma^n_p$ have the Fubini property if, and only, $p=q$.
\end{proposition}

\begin{remark}    \label{R2.6}
Details, explanations and discussions of this substantial assertion may be found in \cite[Section 4, pp.\,34--40]{T01} and 
\cite[Section 3.6, pp.\,103--105]{T20}. The Fubini property admits occasionally to transfer assertions about truncation from one
dimension to higher dimensions. All spaces $B^s_p(\rn)$ in \eqref{2.12} with $s > \sigma^n_p$, including $\Cc^s(\rn)$ with $s>0$
and $H^s_p (\rn)$ with $s>0$, including $W^k_p (\rn)$ have the Fubini property.
\end{remark}

\subsection{Haar bases}   \label{S2.3}
A useful tool to disprove strong truncation and perfect truncation  as described so far in the Introduction  are expansions of elements
belonging to related function spaces in terms of Haar bases. We collect some assertions. We follow \cite[Section 3.5, 
pp.\,98--103]{T20} where it is sufficient for our purpose to concentrate on the $B$--spaces. Let $y\in \real$ and
\begin{\eq}   \label{2.24}
h_M (y) =
\begin{cases}
1 &\text{if $0<y<1/2$},\\
-1 &\text{if $1/2 \le y <1$},\\
0 &\text{if $y \not\in (0,1)$}.
\end{cases}
\end{\eq}
Let $h_F (y) = |h_M (y)|$ be the characteristic function of the unit interval $(0,1)$. Let $n\in \nat$ and let
\begin{\eq}   \label{2.25}
G = (G_1, \ldots, G_n) \in G^0 = \{F,M \}^n
\end{\eq}
which means that $G_r$ is either $F$ or $M$. Let
\begin{\eq}   \label{2.26}
G = (G_1, \ldots, G_n) \in G^* = G^j \in \{F,M \}^{n*}, \qquad j \in \nat,
\end{\eq}
which means that $G_r$ is either $F$ or $M$, where $*$ indicates that at least one of the components of $G$ must be an $M$. Then
\begin{\eq}   \label{2.27}
h^j_{G,m} (x) =  \prod^n_{l=1} h_{G_l} \big( 2^j x_l - m_l \big), \qquad G\in G^j, \quad m \in \zn,
\end{\eq}
$x= (x_1, \ldots, x_n)\in \rn$, where (now) $j\in \no$ is the well--known  Haar basis in \rn. Let $0<p,q \le \infty$ and $s\in \real$.
Then $b^s_{p,q} (\rn)$ is the collection of all sequences
\begin{\eq}   \label{2.28}
\lambda = \big\{ \lambda^{j,G}_m \in \comp: \ j \in \no, \ G \in G^j, \ m \in \zn \big\}
\end{\eq}
such that
\begin{\eq}   \label{2.29}
\| \lambda \, | b^s_{p,q} (\rn) \| = \Big( \sum^\infty_{j=0} 2^{j(s- \frac{n}{p})q} \sum_{G \in G^j} \Big( \sum_{m \in \zn} |\lambda^{j,G}_m|^p \Big)^{q/p} \Big)^{1/q}
\end{\eq}
is finite.

\begin{proposition}   \label{P2.7}
Let $n\in \nat$,
\begin{\eq}   \label{2.30}
0<p,q \le \infty, \qquad \max \Big( n \big(\frac{1}{p}-1), \ \frac{1}{p} - 1 \Big) <s< \min \big( \frac{1}{p}, 1 \big).
\end{\eq}
Let $f \in S'(\rn)$. Then $f\in \Bs (\rn)$ if, and only if, it can be represented by
\begin{\eq}   \label{2.31}
f = \sum_{\substack{j\in \no, G\in G^j,\\ m\in \zn}} \lambda^{j,G}_m \, h^j_{G,m}, \qquad \lambda \in b^s_{p,q} (\rn),
\end{\eq}
the unconditional convergence being in $S' (\rn)$. The representation \eqref{2.31} is unique,
\begin{\eq}   \label{2.32}
\lambda^{j,G}_m =  \lambda^{j,G}_m (f) = 2^{jn} \, \big( f, h^j_{G,m} \big)
\end{\eq}
and
\begin{\eq}   \label{2.33}
I: \quad f \mapsto \big\{ \lambda^{j,G}_m (f) \big\}
\end{\eq}
is an isomorphic map of $\Bs (\rn)$ onto $b^s_{p,q} (\rn)$,
\begin{\eq}   \label{2.34}
\big\| f \, | \Bs (\rn) \big\| \sim \big\| \lambda (f) \, | b^s_{p,q} (\rn) \big\|.
\end{\eq}
\end{proposition}

\begin{remark}   \label{R2.8}
This is essentially the $B$--part of \cite[Theorem 3.18, p.\,99]{T20}. There one finds explanations, discussions, references and a
related (more complicated, but not needed here) counterpart for the $F$--spaces.
\end{remark}

The following observation will be of some use for us later on.

\begin{proposition}    \label{P2.9}
Let $\As (\rn)$ with $A \in \{B,F \}$, $s\in \real$ and $0<p,q \le \infty$ be the spaces according to Definition \ref{D2.1}. Let
$\chi_Q$ be the characteristic function of the cube $Q =(0,1)^n$. Then
\begin{\eq}   \label{2.35}
\chi_Q \in B^{1/p}_{p,q} (\rn) \qquad \text{if, and only if, $q=\infty$}
\end{\eq}
and
\begin{\eq}   \label{2.36}
\chi_Q \in F^{1/p}_{p,q} (\rn) \qquad \text{if, and only if, $p=\infty$}.
\end{\eq}
Furthermore, $\chi_Q \in \As (\rn)$ if, and only if, either $s<1/p$ or as in \eqref{2.35}, \eqref{2.36}.
\end{proposition}

\begin{remark}   \label{R2.10}
This coincides with \cite[Proposition 2.50, p.\,50]{T20}. There one finds also proofs and further discussions.
\end{remark}

\subsection{Faber bases}    \label{S2.4}
In addition to Haar bases we use later on Faber bases for distinguished spaces $\As (\real)$ on \real. Although we deal mainly with
the spaces $\Bs (\real)$ we wish to indicate that the situation for the spaces $\Fs (\real)$, $p \not= q$, including the Sobolev spaces
$H^s_p (\real)$, $p \not= 2$, might be quite different. This suggests to adapt first related assertions from \cite{T10} and \cite{T12}.
But (in contrast to Haar bases) we deal afterwards with some new  aspects which are crucial for our later applications which will be
described in some details.

Let $I =(0,1)$ be the unit interval in $\real$ and let $\chi_{j,m} (x)$, $x\in \real$, be the characteristic functions of the intervals
\begin{\eq}   \label{2.37}
I_{j,m} = \big[2^{-j}m, 2^{-j} (m+1)\big), \qquad j\in \no, \quad \text{and} \quad m=0, \ldots, 2^j -1.
\end{\eq}
Let
\begin{\eq}   \label{2.38}
v_{j,m} (x) =
\begin{cases}
2^{j+1} \big(x- 2^{-j}m \big) &\text{if $2^{-j} m \le x < 2^{-j} m + 2^{-j-1}$,} \\
2^{j+1} \big( 2^{-j} (m+1) -x \big) &\text{if $2^{-j}m + 2^{-j-1} \le x < 2^{-j} (m+1)$,} \\
0 & \text{otherwise},
\end{cases}
\end{\eq}
$x \in \real$, $j\in \no$ and $m=0,\ldots,2^j -1$, be the related well--known Faber hat functions on  $I$. Let
\begin{\eq}   \label{2.39}
\big( \Delta^2_h f \big) (x) = f(x+2h) - 2 f(x+h) +f(x), \qquad x \in \real, \quad h \in \real,
\end{\eq}
be the usual second differences. 

\begin{proposition}    \label{P2.11}
{\em (i)} Let $0<p,q \le \infty$ and
\begin{\eq}   \label{2.40}
\frac{1}{p} <s<1 + \min \big( \frac{1}{p}, 1 \big).
\end{\eq}
Then
\begin{\eq}   \label{2.41}
f \in \Bs (\real) \qquad \text{with} \quad \supp f \subset [0,1] = \ol{I}
\end{\eq}
can be represented as
\begin{\eq}   \label{2.42}
f(x) = - \frac{1}{2} \sum^\infty_{j=0} \sum_{m=0}^{2^j -1} \big( \Delta^2_{2^{-j-1}} f \big) \big( 2^{-j} m \big) v_{j,m} (x),
\qquad x \in \real,
\end{\eq}
with
\begin{\eq}   \label{2.43}
\| f \, | \Bs (\real) \| \sim \| f \, | \Bs (I) \|_d
\end{\eq}
where
\begin{\eq}   \label{2.44}
\| f \, | \Bs (I) \|_d = \bigg( \sum^\infty_{j=0} 2^{j(s- \frac{1}{p})q} \Big( \sum^{2^j -1}_{m=0} \big| \big(\Delta^2_{2^{-j-1}} f
\big) (2^{-j} m ) \big|^p \Big)^{q/p} \bigg)^{1/q}
\end{\eq}
$($equivalent quasi--norms$)$, usual modification if $\max(p,q) = \infty$.
\cm
{\em (ii)} Let
\begin{\eq}   \label{2.45}
\begin{cases}
0<p<\infty, 0<q <\infty, & \max \big( \frac{1}{p}, \frac{1}{q}, 1 \big) <s< 1 + \min \big( \frac{1}{p}, \frac{1}{q}, 1 \big), \\
1<p<\infty, 1<q<\infty, &s=1, \\
1<p<\infty, 1<q \le \infty, &\max \big( \frac{1}{p}, \frac{1}{q} \big) <s<1.
\end{cases}
\end{\eq}
Then
\begin{\eq}   \label{2.46}
f \in \Fs (\real) \qquad \text{with} \quad \supp f \subset [0,1] = \ol{I}
\end{\eq}
can be represented by \eqref{2.42} with
\begin{\eq}   \label{2.47}
\|f \, | \Fs (\real) \| \sim \| f \, | \Fs (I) \|_d
\end{\eq}
where
\begin{\eq}   \label{2.48}
\| f \, | \Fs (I) \|_d = \bigg\| \Big( \sum^\infty_{j=0} \sum^{2^j -1}_{m=0} 
2^{jsq} \big| \big(\Delta^2_{2^{-j-1}} f
\big) (2^{-j} m ) \big|^q \chi_{j,m} (\cdot) \Big)^{1/q} | L_p (I) \bigg\|
\end{\eq}
$($equivalent quasi--norms$)$, usual modification if $q=\infty$.
\end{proposition}

\begin{remark}   \label{R2.12}
This assertion is covered  by \cite[Theorem 3.1, pp.\,126--127]{T10} and \cite[Theorem 3.6, pp.\,49--50]{T12}. There one finds also
explanations and proofs. This means 
\begin{\eq}   \label{2.49}
\begin{cases}
2\le p <\infty, & \frac{1}{2} <s< 1 + \frac{1}{p}, \\
p<2, & \frac{1}{p} <s< \frac{3}{2},
\end{cases}
\end{\eq}
for the (extended) Sobolev spaces $H^s_p(\real) = F^s_{p,2} (\real)$ covering in particular corresponding Sobolev spaces with $1<p<\infty$ according to \eqref{2.14}. At the first glance the conditions 
\eqref{2.45} and also \eqref{2.49} may look a little bit artificial (depending on the method). But this is not the case. They are 
natural. This follows from  the discussions about Haar bases in \cite[Section 3.5, pp.\,98--103]{T20}, based on related references,
and the observation that restrictions for the smoothness $s$ for Haar bases and Faber bases differ by 1.
\end{remark}

One may ask for a counterpart of Proposition \ref{P2.9} in one dimension for the basic Faber function
\begin{\eq}   \label{2.50}
v(x) = \max( 1-|x|,0), \qquad x\in \real,
\end{\eq}
underlying the above considerations, based on \cite[Section 3.1.2, pp.\,126--129]{T10} and \cite[Section 3.2, pp.\,44--51]{T12}.

\begin{proposition}   \label{P2.13}
Let $\As (\real)$ with $A \in \{B,F \}$, $s\in \real$ and $0<p,q \le \infty$ be the spaces according to Definition \ref{D2.1} with 
$n=1$. Then
\begin{\eq}  \label{2.51}
v \in B^{1+ \frac{1}{p}}_{p,q} (\real) \qquad \text{if, and only if, $q=\infty$},
\end{\eq}
and
\begin{\eq}   \label{2.52}
v\in F^{1+ \frac{1}{p}}_{p,q} (\real) \qquad \text{if, and only if, $p=\infty$}.
\end{\eq}
Furthermore, $v\in \As (\real)$ if, and only if, either $s<1 + \frac{1}{p}$ or as in \eqref{2.51}, \eqref{2.52}.
\end{proposition}

\begin{remark}   \label{R2.14}
Compared with Proposition \ref{P2.9} one has to lift the smoothness $s$ by $1$. This is plausible. But the detailed justification
requires some care, \cite[Section 5.8]{Tri21}.
\end{remark}

We wish to use the above considerations to discuss strong and perfect truncations in one dimension as described so far in the 
Introduction. For this purpose one has to complement and adapt what has been recalled above. We rely again on the related parts in
\cite{T10} in the version as developed  in \cite{T12}, where we extended the Faber expansions from $I = (0,1)$ to \real. Let again
$0<p,q \le \infty$ and
\begin{\eq}   \label{2.53}
\frac{1}{p} <s< 1 + \min \big( \frac{1}{p}, 1 \big).
\end{\eq}
Let $f\in \Bs (\real)$ with $f(k) =0$ if $k\in \ganz$ and compact support. Then it follows from \cite[Theorem 3.6, pp.\,49--50]{T12}
that
\begin{\eq}   \label{2.54}
f(x) = - \frac{1}{2} \sum^\infty_{j=0} \sum_{m \in \ganz} \big( \Delta^2_{2^{-j-1}} f \big) \big( 2^{-j} m \big) v_{j,m} (x),
\qquad x \in \real,
\end{\eq}
extending \eqref{2.38} from $I = (0,1)$ to \real, with
\begin{\eq}   \label{2.55}
\| f \, | \Bs (\real) \| \sim \| f \, | \Bs (\real) \|_d
\end{\eq}
where
\begin{\eq}   \label{2.56}
\| f \, | \Bs (\real) \|_d = \bigg( \sum^\infty_{j=0} 2^{j(s- \frac{1}{p})q} \Big(
\sum_{m \in \ganz} \big| \big(\Delta^2_{2^{-j-1}} f \big) (2^{-j} m ) \big|^p \Big)^{q/p} \bigg)^{1/q}
\end{\eq}
(with the usual modification  if $\max(p,q) = \infty$). If $p=q$ then one has even
\begin{\eq}   \label{2.57}
\| f \, | B^s_p (\real) \|^p_d = \sum_{k\in \ganz} \|f \, | B^s_p (k + I)\|^p_d
\end{\eq}
in obvious extension of \eqref{2.44}. But this will not be needed.
Let, in addition, $f$ be real. Then it follows from the already well established  truncation property recalled in Theorem 
\ref{T3.12}(i) below that
\begin{\eq}   \label{2.58}
\| f^+ \, | \Bs (\real) \|  + \| f^- \, | \Bs (\real) \| \le c \, \| f \, | \Bs (\real) \|
\end{\eq}
where
\begin{\eq}   \label{2.59}
f^+ (x) = \max \big( f(x), 0 \big) \quad \text{and} \quad f^- (x) = \min \big( f(x), 0 \big), \qquad x \in \real.
\end{\eq}
Recall that \eqref{2.53} ensures that $f$ is continuous. Furthermore $f(x) = f^+ (x) + f^- (x)$ and $|f(x)|
= f^+ (x) - f^- (x)$. Let us suppose in addition to the above assumptions including $f(k) =0$ if $k\in \ganz$ that the real function
$f(x)$ does not
change sign in any interval $k + I = (k,k+1)$, $k \in \ganz$. Then it follows from the above considerations that
\begin{\eq}   \label{2.60}
\| f \, | \Bs (\real) \|_d \sim \| f^+ \, | \Bs (\real) \|_d + \|f^- \,| \Bs (\real) \|_d \sim
\big\| |f| \, | \Bs (\real) \big\|_d.
\end{\eq}
The question arises to which extent \eqref{2.60} rewritten as
\begin{\eq}   \label{2.61}
\| f \, | \Bs (\real) \| \sim \| f^+ \, | \Bs (\real) \| + \|f^- \,| \Bs (\real) \| \sim
\big\| |f| \, | \Bs (\real) \big\|
\end{\eq}
remains valid for arbitrary
\begin{\eq}   \label{2.62}
f\in\Bs (\real) \ \text{real}, \quad 0<p,q \le \infty, \quad \frac{1}{p} <s<1 + \min \big( \frac{1}{p}, 1 \big),
\end{\eq}
where the related equivalence constants are independent of $f$. But this is the case and can be justified as follows. Let
\begin{\eq}   \label{2.63}
f \in \Bs (\real), \qquad \supp f \subset [0, \lambda] = I_\lambda, \qquad 0<\lambda \le 1,
\end{\eq}
with $p,q,s$ as in \eqref{2.62}. Then one has by the local homogeneity as described in \cite[Theorem 3.9, pp.\,91--92]{T20} and the
references given there that 
\begin{\eq}   \label{2.64}
\| f(\lambda \cdot) \, | \Bs (\real) \| \sim \lambda^{s- \frac{1}{p}} \, \| f \, | \Bs (\real) \|,
\end{\eq}
where the equivalence constants are independent of $f$ and $\lambda$ satisfying \eqref{2.63}. Then can apply Proposition \ref{P2.11}(i)
to $f(\lambda \cdot)$. Re--transformation gives a Faber representation if $\lambda = 2^{-J}$, $J\in \nat$, consisting exclusively
of terms $v_{j,m}$ with $j \ge J$.
Otherwise one has immaterial perturbations which do not influence the arguments. This shows how to deal with small 
intervals between consecutive zeros of $f(x)$. If the related intervals are large then one can rely in addition on \cite[Theorem 3.6,
pp.\,49--50]{T12}. Finally one has to clip together these intervals. To avoid awkward situations (accumulation points of zeros of 
$f(x)$) one may assume first $f=g \psi$, where $g$ is an analytic function (going back to the functions in the original Definition
\ref{D2.1}(i)) and $\psi$ is a smooth cut--off function. Afterwards one can rely on completion and Fatou arguments.

\subsection{Oscillations}    \label{S2.5}
A further useful instrument to deal with strong and perfect  truncation is the characterization of some function spaces $\As (\rn)$ in
terms of oscillations. Let $1 \le r \le \infty$. Then
\begin{\eq}   \label{2.65}
\sigma^{r,n}_p = n \Big( \max \big( \frac{1}{p}, \frac{1}{r} \big) -\frac{1}{r} \Big), 
\qquad \sigma_{p,q}^{r,n} = n \Big( \max \big( \frac{1}{p}, \frac{1}{q}, \frac{1}{r} \big) - \frac{1}{r} \Big)
\end{\eq}
generalizes $\sigma^{1,n}_p = \sigma^{n}_p$ and $\sigma^{1,n}_{p,q} = \sigma^{n}_{p,q}$ according to \eqref{2.17}, where again
$0<p,q \le \infty$ and $n\in \nat$. Recall that $L^{\loc}_u (\rn)$ consists  of all Lebesgue measurable (complex--valued) functions
$f$ on $\rn$ such that $f\in L_u (\Om)$ for any bounded domain $\Omega$ in $\rn$ where $0<u \le \infty$. 
Let $f\in L^{\loc}_u (\rn)$ where $0<u \le \infty$ and $M\in \no$. Then
\begin{\eq}   \label{2.66}
\osc^M_u f(x,t) = \inf \Big( t^{-n} \int_{\{y: |x-y| <t \}} | f(y) - P(y) |^{u} \, \di y \Big)^{1/u}, \qquad x\in \rn, \quad t>0,
\end{\eq}
where the infimum is taken over all polynomials $P$ of degree less than or equal to $M$, is called the (local) oscillation  of $f$.
We dealt in detail in \cite[Section 3.4, pp.\,179--186]{T92} with oscillations. They can be used to characterize some spaces $\As (\rn)
$. We follow here \cite[Section 1.7.3, pp.\,51--52]{T92} which, in turn, is based on \cite{Tri89}. There one finds also related 
references. If $s\in \real$ then $[s] \in \ganz$ is the largest integer smaller than or equal to $s$.

\begin{proposition}   \label{P2.15}
{\em (i)} Let $0<p\le \infty$, $0<q \le \infty$, $1\le r \le \infty$ and $s> \sigma^{r, n}_p$. Let $0<u \le r$ and $[s] \le M \in 
\no$. Then $\Bs (\rn)$ is the collection of all $f\in L_{\max(p,r)} (\rn)$ such that
\begin{\eq}   \label{2.67}
\| f \, | L_p (\rn) \| + \Big( \int^1_0 t^{-sq} \big\| \osc^M_u f (\cdot, t) \, | L_p (\rn) \big\|^q \, \frac{\di t}{t} \Big)^{1/q}
\end{\eq}
is finite $($equivalent  quasi--norms$)$, usual modification if $q=\infty$.
\cm
{\em (ii)} Let $0<p<\infty$, $0<q \le \infty$, $1\le r \le \infty$ and $s > \sigma^{r,n}_{p,q}$. Let $0<u \le r$ and $[s] \le M \in
\no$. Then $\Fs (\rn)$ is the collection of all $f\in L_{\max(p,r)} (\rn)$ such that
\begin{\eq}    \label{2.68}
\| f \, | L_p (\rn) \| + \Big\| \Big( \int^1_0 t^{-sq} \osc^M_u f (\cdot, t)^q  \, \frac{\di t}{t} \Big)^{1/q} \, | L_p (\rn) \Big\|
\end{\eq}
is finite $($equivalent quasi--norms$)$, usual modification if $q=\infty$.
\end{proposition}

\begin{remark}   \label{R2.16}
A detailed proof may be found in \cite[Section 3.5.1, pp.\,186--192]{T92} where in the corresponding theorem the misprint $M > [s]$
in connection with part (ii) of the above proposition must be corrected by $M \ge [s]$ in agreement as already correctly stated in
\cite[Theorem 1.7.3, p.\,51]{T92}.
\end{remark}

Of interest for us is the case $r= u= \infty$. Then the above proposition requires $s> n/p$ ensuring in particular
\begin{\eq}    \label{2.69}
 \As (\rn) \hra C(\rn).
\end{\eq}
Let $f \in C(\rn)$ be real. Then (in modification of \eqref{2.66} with $u=\infty$)
\begin{\eq}    \label{2.70}
[f]_{j,m} = \sup_{y\in Q_{j,m}} f(y) - \inf_{y \in Q_{j,m}} f(y) = 2 \, \osc^0_\infty f(Q_{j,m}), \qquad j \in \no, \quad m\in \zn,
\end{\eq}
with $Q_{j,m} = 2^{-j} m + 2^{-j} (0,1)^n$ as in \eqref{2.6}.
This can be applied to the above Proposition \ref{2.15} if, in addition, $0<s<1$. Then the second terms in \eqref{2.67} and 
\eqref{2.68} can be discretized. Let
\begin{\eq}   \label{2.71}
\big\| [f] \, | b^s_{p,q} (\rn) \big\| = \bigg( \sum^\infty_{j=0} 2^{j(s- \frac{n}{p})q} \Big( \sum_{m\in \zn} \big| [f]_{j,m} \big|^p
\Big)^{q/p} \bigg)^{1/q}
\end{\eq}
and
\begin{\eq}   \label{2.72}
\big\| [f] \, | f^s_{p,q} (\rn) \big\| = \Big\| \Big( \sum^\infty_{j=0} \sum_{m \in \zn} 2^{jsq} \big| [f]_{j,m} \chi_{j,m} (\cdot)
\big|^q \Big)^{1/q} \, | L_p (\rn) \Big\|
\end{\eq}
where
\begin{\eq}   \label{2.73}
[f] = \big\{ [f]_{j,m} : \ j \in \no, \ m \in \zn \big\}.
\end{\eq}
Here again $\chi_{j,m}$ is the characteristic function of $Q_{j,m} = 2^{-j}m + 2^{-j} (0,1)^n$. If $f\in C(\rn)$ is complex then
\begin{\eq}   \label{2.74}
\big\| [f] \, | a^s_{p,q} (\rn) \big\| = \big\| [\Rea f] \, | a^s_{p,q} (\rn) \big\| + \big\| [\Ima f] \, | a^s_{p,q} (\rn) \big\|.
\end{\eq}
We fix the outcome.  

\begin{theorem}   \label{T2.17}
{\em (i)} Let $0<p,q \le \infty$ and $n/p <s<1$. Then $\Bs (\rn)$ is the collection of all $f\in C(\rn)$ such that
\begin{\eq}   \label{2.75}
\|f \, | L_p (\rn) \| + \big\| [f] \, | b^s_{p,q} (\rn) \big\|
\end{\eq}
is finite $($equivalent quasi--norms$)$.
\cm
{\em (ii)} Let $0<p<\infty$, $0<q \le \infty$ and $n \max \big( \frac{1}{p}, \frac{1}{q} \big) <s<1$. Then $\Fs (\rn)$ is the 
collection of all $f\in C(\rn)$ such that
\begin{\eq}   \label{2.76}
\|f \, | L_p (\rn) \| + \big\| [f] \, | f^s_{p,q} (\rn) \big\|
\end{\eq}
is finite $($equivalent norms$)$.
\end{theorem}

\begin{remark}    \label{R2.18}
This follows from Proposition \ref{P2.15} and the above considerations. The discretization of $t$ with $0<t<1$ by $\{ 2^{-j}: j\in \no\}$ had already been done in \cite[(6), (7), p.\,187]{T92}. This can be easily complemented by the indicated discretization of $\rn$ in $\{ Q_{j,m}: m\in \zn \}$, $j\in \no$. 
\end{remark}

Although the above theorem looks very special, it might be of some self--contained interest. In any case it paves the way to say 
something about diverse truncation properties as mentioned so far in the Introduction. The continuous
version of \eqref{2.70} inserted in \eqref{2.67} with $n=1$, $1<p=q <\infty$, $1/p <s<1$, $r=u= \infty$, coincides essentially with
\cite[(1.6), (1.7)]{MiVS19}, where one finds a related  short proof. 

\section{Truncations}   \label{S3}
\subsection{Definitions, basic assertions and examples}    \label{S3.1}
Let again $C(\rn)$ be the naturally normed space of bounded continuous functions on \rn. Let
\begin{\eq}   \label{3.1}
\sigma^n_p = n \Big( \max \big( \frac{1}{p}, 1 \big) - 1 \Big) \quad \text{and} \quad \sigma^n_{p,q} = n \Big( \max \big(\frac{1}{p}, \frac{1}{q}, 1 \big) - 1 \Big),
\end{\eq}
where $n\in \nat$ and $0<p,q \le \infty$ be as in \eqref{2.17}. Truncation applies to spaces $\As (\rn)$ consisting entirely of 
regular distributions. This means \eqref{2.18} specified by \eqref{2.19} and some limiting cases with $s= \sigma^n_p$. Let $\As (\rn)$
be the spaces as introduced in Definition \ref{D2.1}, $n\in \nat$.

\begin{definition}   \label{D3.1}
{\em (i)} Let $\As (\rn) \subset L^{\loc}_1 (\rn)$. Then $\As (\rn)$ is said to have the truncation property if $|f| \in \As (\rn)$ for
any real $f\in \As (\rn)$ and if there is a constant $c>0$ such that
\begin{\eq}    \label{3.2}
\big\| \, |f| \, \big| \As (\rn) \big\| \le c \, \| f \, | \As (\rn) \|, \qquad f \in \As (\rn) \ \text{real}.
\end{\eq}
{\em (ii)} Let $\As (\rn) \subset L^{\loc}_1 (\rn)$. Then $\As (\rn)$ is said to have the strong truncation property if $|f| \in \As
(\rn)$ for any real $f\in \As (\rn)$ and
\begin{\eq}    \label{3.3}
\big\| \, |f| \, \big| \As (\rn) \big\| \sim \| f \, | \As (\rn) \|, \qquad f \in \As (\rn) \ \text{real},
\end{\eq}
$($equivalences$)$.
\cm
{\em (iii)} Let $\As (\rn) \hra C(\rn)$ $($continuous embedding$)$. 
Then $\As (\rn)$ is said to have the perfect truncation property when any real $f\in C(\rn)$ belongs to $\As (\rn)$ if, and only if, $|f|$ belongs to $\As (\rn)$ and
\begin{\eq}    \label{3.4}
\big\| \, |f| \, \big| \As (\rn) \big\| \sim \| f \, | \As (\rn) \|,
\end{\eq}
$($equivalences$)$. 
\end{definition}

\begin{remark}   \label{R3.2}
The restriction $f\in L^{\loc}_1 (\rn)$ in parts (i) and (ii) of the above definition is natural. The assumption $\As (\rn) \hra C(\rn)
$ and $f \in C(\rn)$ real in part (iii) seems to be at least reasonable. It underlies the arguments in respective assertions. But one
may ask whether $C(\rn)$ can be replaced by weaker requests, for example $L^{\loc}_1 (\rn)$. But this may cause some extra problems 
and will not be considered here. It may happen that $|f|$ is much smoother than $f$. This can be illustrated by the examples below.
Instead  of $|f|$ for $f$ real one may work with $f^+$ and $f^-$,
\begin{\eq}   \label{3.5}
f^+ (x) = \max (f(x), 0), \qquad f^- (x) = \min (f(x), 0), \qquad \text{$x\in \rn$ a.e.,}
\end{\eq}
as we already did in \eqref{2.58}, \eqref{2.59}, based on
\begin{\eq}    \label{3.6}
f(x) = f^+ (x) + f^- (x) \quad \text{and} \quad |f(x)| = f^+ (x) - f^- (x).
\end{\eq}
But we stick preferably on $|f|$. More detailed discussions may be found in \cite[pp.\,357--362]{T01}.
\end{remark}

We considered in \cite[Chapter 25]{T01} truncations with the following outcome. Let $\sigma^n_p$ and $\sigma^n_{p,q}$ be as in 
\eqref{3.1}.

\begin{theorem}   \label{T3.3}
Let $n\in \nat$.
\cm
{\em (i)} The spaces  $\Bs (\rn)$ have for all $q$ with $0<q \le \infty$ the truncation property, if, and only if,
\begin{\eq}   \label{3.7}
0<p \le \infty \qquad \text{ and} \qquad \sigma^{n}_p <s< 1 + \frac{1}{p}.
\end{\eq}
{\em (ii)} The spaces $\Fs (\rn)$ have the truncation property if
\begin{\eq}   \label{3.8}
0<p<\infty, \quad 0<q \le \infty \quad \text{and} \quad \sigma^n_{p,q} <s<1 + \frac{1}{p}
\end{\eq}
$( s \not= \frac{1}{p}$ for $0<p \le 1)$.
\end{theorem}

\begin{remark}   \label{R3.4}
This is covered by \cite[Theorem 25.8, p.\,364]{T01}. The proof is rather long, some 10 pages, relying on atomic decompositions, where 
especially the cases $s>1$, $s>2$ and $0<p <1$ required some specific arguments. We could not remove
the unnatural restriction  $s \not= 1/p$ if $0<p \le 1$. It depends surely on the method. An extension of this assertion to so--called
tempered homogeneous spaces may be found in \cite[Section 3.11]{T15}.
\end{remark}

We discuss some special cases and aspects. Further examples may be found in  \cite[Section 25.1, pp.\,357--361]{T01}.

\begin{example}  \label{E3.5}
Obviously,
\begin{\eq}   \label{3.8a}
\big\| \, |f| \, | L_p (\rn) \big\| = \| f \, | L_p (\rn) \|, \qquad 1 \le p \le \infty.
\end{\eq}
It is less obvious that this assertion can be extended to the classical Sobolev spaces $W^1_p (\rn)$, normed by \eqref{2.16}, $1 \le p
<\infty$ (including temporarily $p=1$). Let $f \in W^1_p (\rn)$ be real. Then
\begin{\eq}    \label{3.9}
\big\| \, |f| \, | W^1_p (\rn) \big\| = \| f \, | W^1_p (\rn) \|.
\end{\eq}
This is covered by \cite[Corollary 2.1.8, p.\,47]{Zie89} and also \cite[Lemma 7.6, p.\,145]{GiT77}.
\end{example}

\begin{example}   \label{E3.6}
We discuss a further simple example. The H\"{o}lder spaces $\Cc^s (\rn)$ with $0<s<1$ according to \eqref{2.13} can be equivalently
normed by
\begin{\eq}   \label{3.10}
\| f \, | \Cc^s (\rn) \| = \| f\, | L_\infty (\rn) \| + \sup_{0<|x-y|\le 1} \frac{|f(x) - f(y)|}{|x-y|^s}.
\end{\eq}
Then $|f| \in \Cc^s (\rn)$ if $f\in \Cc^s(\rn)$, $0<s<1$, and
\begin{\eq}   \label{3.11}
\big\| \, |f| \, |\Cc^s(\rn) \big\| \le \| f \, | \Cc^s (\rn) \|.
\end{\eq}
There is a converse. Let $f\in C(\rn)$ be real and $|f| \in \Cc^s (\rn)$. Let $f(x) >0$, $f(y) <0$ and $f(z) =0$ for some $z= (1-t)x
+ty$ with $0<t<1$. Then it follows from
\begin{\eq}   \label{3.12}
\frac{ f(x) - f(y) }{|x-y|^s} = \frac{ |f(x)| + |f(y)| }{|x-y|^s} \le \frac{ |f(x)| - |f(z)| }{|x-z|^s}+
\frac{ |f(y)| - |f(z)| }{|y-z|^s}
\end{\eq}
that $f\in \Cc^s (\rn)$ and, together with \eqref{3.11},
\begin{\eq}   \label{3.13}
\| f \, | \Cc^s (\rn) \| \sim \big\| \, |f| \, | \Cc^s (\rn) \big\|, \qquad f \in \Cc^s(\rn) \ \text{real}.
\end{\eq}
In other words, the spaces $\Cc^s (\rn)$ with $0<s<1$ have the perfect truncation property according to Definition \ref{D3.1}(iii).
\end{example}

\begin{example}    \label{E3.7}
On the other hand we wish to illuminate the significant difference between  truncation, strong truncation and perfect truncation as
introduced in Definition \ref{D3.1}. Let $n\in \nat$,
\begin{\eq}   \label{3.14}
0<p,q \le \infty \quad \text{and} \quad \max \Big( n \big( \frac{1}{p} -1 \big), \frac{1}{p} -1 \Big) <s< \min \big( \frac{1}{p}, 1 \big).
\end{\eq}
Then it follows from the Haar expansion according to Proposition \ref{P2.7} that 
\begin{\eq}   \label{3.15}
f_j (x) = \sum_{m_l =0,..,2^j -1} h^j_{G,m} (x) \in B^s_{p,q} (\rn), \qquad j\in \nat,
\end{\eq}
$G= (M,\ldots, M)$,  and
\begin{\eq}   \label{3.16}
\| f_j \, | B^s_{p,q} (\rn) \| \sim 2^{js}, \qquad j \in \nat,
\end{\eq}
where the equivalence constants are independent of $j$. Furthermore, $|f_j| = \chi_Q$ is the characteristic function of $Q= (0,1)^n$. In  particular 
\begin{\eq}   \label{3.17}
\big\| \, |f_j | \, | B^s_{p,q} (\rn) \big\| \sim 1, \qquad j\in \nat.
\end{\eq}
It follows from
\begin{\eq}   \label{3.18}
B^s_{p,\min(p,q)} (\rn) \hra \Fs (\rn) \hra B^s_{p, \max(p,q)} (\rn)
\end{\eq}
according to \cite[Theorem 2.9, p.\,26]{T20} (including $F^s_{\infty,q} (\rn)$) that one can replace $B$ in \eqref{3.16} and 
\eqref{3.17} by $A$ with $A\in \{B,F \}$. If
\begin{\eq}   \label{3.19}
n \Big( \max \big( \frac{1}{p}, 1 \big) - 1 \Big) = \sigma^{(n)}_p <s< \min \big(\frac{1}{p}, 1 \big), \qquad 0<p<\infty,
\end{\eq}
$0<q \le \infty$. Then one has
\begin{\eq}   \label{3.20}
\big\| \, |f_j | \, | \As (\rn) \big\| \sim 1 <2^{js} \sim \| f_j \, | \As (\rn) \|, \qquad j \in \nat.
\end{\eq}
This shows that those spaces covered also by Theorem \ref{T3.3}, this means $\Bs (\rn)$ with \eqref{3.14} and its $F$--counterpart
based on \eqref{3.8} and $s< \min(\frac{1}{p},1)$ have the truncation property but not the strong truncation property.
\end{example}

The above considerations can be extended and formalized as follows. Let again $\chi_Q$ be the characteristic function of the cube $Q=
(0,1)^n$ and let $\sigma^n_p$ be as in \eqref{3.1}.

\begin{proposition}    \label{P3.8}
Let $0<p<\infty$, $0<q \le \infty$ and $s > \sigma^{n}_p$. If $\chi_Q \in \As (\rn)$ then $\As (\rn)$ has not the strong truncation
property.
\end{proposition}

\begin{remark}   \label{R3.9}
This follows essentially from the considerations  in Example \ref{E3.7}, the embedding
\begin{\eq}    \label{3.21}
B^{s_1}_{p_1,q} (\rn) \hra B^{s_0}_{p_0,q} (\rn), \qquad s_1 - \frac{n}{p_1} = s_0 - \frac{n}{p_0},
\end{\eq}
where $0<p_1 \le p_0 <\infty$, $0<q \le \infty$ and $s_1 > \sigma^n_{p_1}$, and Proposition \ref{P2.9}. This proposition justifies
also the following assertion.
\end{remark}

\begin{corollary}   \label{C3.10}
Let $n\in \nat$. The spaces
\begin{\eq}    \label{3.22}
\As (\rn), \qquad 0<p<\infty, \quad 0<q \le \infty, \quad \sigma^{n}_p <s< \frac{1}{p},
\end{\eq}
and $B^{1/p}_{p,\infty} (\rn)$, $\frac{n-1}{n} <p<\infty$, do not have the strong truncation property.
\end{corollary}

\begin{remark}  \label{R3.11}
According to \cite[Proposition 4.2]{VScha20} the spaces $B^{1/p}_p (\rn) = B^{1/p}_{p,p} (\rn)$, $1<p<\infty$,  have also not the strong truncation property. These spaces are not covered by the above corollary.
What about the spaces $B^{1/p}_{p,q} (\rn)$ with $p\not=q$? We return below to this point in Theorem \ref{T3.12}(iii) with the
expected outcome that, at least in one dimension, all related spaces $B^{1/p}_{p,q} (\real)$, $0<q \le \infty$, do not have the strong
truncation property.
\end{remark}

\subsection{Truncations in one dimension}   \label{S3.2}
Truncations have been introduced in Definition \ref{D3.1}. If the space $\As (\rn)$ has the Fubini property according to Proposition
\ref{P2.5} then one can transfer strong truncation  from one dimension to higher dimensions. But this is not the case for perfect
truncations. This suggests to deal separately  with truncations in one dimension and in higher dimensions. Furthermore assertions for
$F$--spaces are less complete than for their $B$--counterparts, occasionally even based on different arguments. This may justify to 
concentrate first on $B$--spaces  and then on $F$--spaces. We rely on the tools as provided in Section \ref{S2} and hint how and where
to use them. But for further details we refer the reader again to the Informal Notes \cite{Tri21}.

\begin{theorem}   \label{T3.12}
{\em (i)} The spaces $\Bs (\real)$ with 
\begin{\eq}   \label{3.23}
 0<p,q \le \infty \quad \text{and} \quad \max \Big( \frac{1}{p}, 1 \Big) -1 <s< 1 + \frac{1}{p}
\end{\eq}
have the truncation property.
\cm
{\em (ii)} The spaces $\Bs (\real)$ with 
\begin{\eq}   \label{3.24}
0<p,q \le \infty \quad \text{and} \quad \frac{1}{p} <s<1 + \min \Big( \frac{1}{p}, 1 \Big)
\end{\eq}
have the perfect truncation property.
\cm
{\em (iii)} 
The spaces $\Bs (\real)$ with 
\begin{\eq}   \label{3.25}
 0<p<\infty, \ 0<q \le \infty \quad \text{and} \quad \max \Big( \frac{1}{p}, 1 \Big) -1 <s \le
\frac{1}{p}
\end{\eq}
do not have the strong truncation property.
\end{theorem}

\begin{remark}   \label{R3.13}
Part (i) follows from Theorem \ref{T3.3}(i) with $n=1$. The proof of part (ii) can be based on the Faber bases as considered in Section
\ref{S2.4} and is essentially covered by \eqref{2.60} and the discussions afterwards. Compared with  Corollary \ref{C3.10} it remains
to justify that the spaces $B^{1/p}_{p,q} (\real)$ with $0<q<\infty$ do not have the strong truncation property. But this is covered 
by \cite[Proposition 5.36, p.\,84]{Tri21} where we relied on characterizations of $B^{1/p}_{p,q} (\real)$ in terms of differences
applied to some special functions. We do not go into the details.
\end{remark}

The $F$--counterpart of the above theorem is less complete.

\begin{theorem}   \label{T3.14}
{\em (i)} The spaces $\Fs (\real)$ with
\begin{\eq}   \label{3.26}
0<p<\infty, \quad 0<q \le \infty \quad \text{and} \quad \max \big( \frac{1}{p}, \frac{1}{q},1 \big) -1 <s< 1 + \frac{1}{p}
\end{\eq}
$( s \not= \frac{1}{p}$ for $0<p \le 1)$ have the truncation property.
\cm
{\em (ii)} The spaces $\Fs (\real)$ with
\begin{\eq}   \label{3.27}
1<p<\infty, \quad 1<q \le \infty \quad \text{and} \quad \max \big( \frac{1}{p}, \frac{1}{q} \big) <s<1
\end{\eq}
have the perfect truncation property.
\cm
{\em (iii)} The spaces $\Fs (\real)$ with
\begin{\eq}   \label{3.28}
0<p<\infty, \quad 0<q \le \infty \quad \text{and} \quad \max \big( \frac{1}{p}, 1 \big) -1 <s< \frac{1}{p}
\end{\eq}
do not have the strong truncation property.
\end{theorem}

\begin{remark}    \label{R3.15}
Part (i) follows from Theorem \ref{T3.3}(ii) with $n=1$. In contrast to the $B$--spaces related Faber expansions  for the spaces $\Fs(
\real)$ as described in Proposition \ref{P2.11}(ii) cannot be used to prove part (ii). But rescue comes from the characterization of
$\Fs (\real)$ with \eqref{3.27} in terms of oscillations according to Theorem \ref{T2.17}(ii) and
\begin{\eq}   \label{3.29}
[f]_{j,m} \sim [|f|]_{j,m}, \qquad j\in \no, \quad m\in \ganz,
\end{\eq}
based on \eqref{2.70}. This applies also to $B$--spaces, but there we have already in Theorem \ref{T3.12}(ii) a better assertion. 
Finally part (iii) is a special case of Corollary \ref{C3.10} with $n=1$.
\end{remark}

\subsection{Truncations in higher dimensions}    \label{S3.3}
One can transfer truncation and strong truncation as introduced in Definition \ref{D3.1} from one dimension to higher dimensions if the
underlying spaces have the Fubini property according to  Proposition \ref{P2.5}. This applies to the $F$--spaces  restricted by 
\eqref{2.23} and does not apply to the $B$--spaces $\Bs (\rn)$ if $p \not= q$. Perfect truncations cannot be shifted from one dimension
to higher dimensions in this way, because \eqref{2.22} requires  that one already knows that $f\in \As (\rn)$. Furthermore we 
restricted perfect truncation in Definition \ref{D3.1}(iii) to spaces which are continuously embedded in $C(\rn)$. If the Fubini
property cannot applied  then one has to rely on direct arguments based on the instruments as provided in Section \ref{S2} with the
following outcome, separating again $B$--spaces and $F$--spaces. Let, as before, for $n\in \nat$,
\begin{\eq}   \label{3.30}
\sigma^n_p = n \Big( \max \big( \frac{1}{p}, 1 \big) - 1 \Big) \quad \text{and} \quad 
\sigma^n_{p,q} = n \Big( \max \big(\frac{1}{p}, \frac{1}{q},
1 \big) - 1 \Big), 
\end{\eq}
and $B^s_p (\rn) = B^s_{p,p} (\rn)$, $s\in \real$, where $0<p,q \le \infty$.

\begin{theorem}   \label{T3.16}
Let $n \in \nat$.
\cm
{\em (i)} The spaces $\Bs (\rn)$ with
\begin{\eq}   \label{3.31}
 0<p,q \le \infty \quad \text{and} \quad \sigma^{n}_p <s< 1 + \frac{1}{p}
\end{\eq}
have the truncation property.
\cm
{\em (ii)} The spaces $\Bs (\rn)$ with
\begin{\eq}   \label{3.32}
 0<p,q \le \infty \quad \text{and} \quad \frac{n}{p} <s<1
\end{\eq}
have the perfect truncation property.
\cm
{\em (iii)} The spaces $B^s_p (\rn)$ with
\begin{\eq}   \label{3.33}
 0<p \le \infty \quad \text{and} \quad \max \Big(\frac{1}{p}, \sigma^{n}_p \Big)
 <s< 1 + \min \Big( \frac{1}{p}, 1 \Big)
\end{\eq}
have the strong truncation property.
\cm
{\em (iv)} The spaces $\Bs (\rn)$ with
\begin{\eq}   \label{3.34}
 0<p<\infty, \quad 0<q \le \infty \quad \text{and} \quad \sigma^{n}_p <s \le \frac{1}{p},
\end{\eq}
do not have the strong truncation property.
\end{theorem}

\begin{remark}    \label{R3.17}
This extends Theorem \ref{T3.12} naturally from $\real$ to \rn, $n\in \nat$. Part (i) is again covered by Theorem \ref{T3.3}(i). In  
contrast to the proof of Theorem \ref{T3.12}(ii) where we used Faber expansions in $\real$ (having no higher--dimensional counterpart)
one has now to rely on the characterization of $\Bs (\rn)$ in terms of oscillations according to Theorem \ref{T2.17}(i) and again 
\eqref{3.29} based on \eqref{2.70}. Part (iii) for $B^s_p (\rn) = B^s_{p,p} (\rn) = F^s_{p,p} (\rn)$ with $p<\infty$ follows from
Theorem \ref{T3.12}(ii) and the Fubini property according to Proposition \ref{P2.5}(i). The case $\Cc^s (\rn) = B^s_\infty (\rn)$ with
$0<s<1$ is already covered by part (ii) and also by \eqref{3.13}. Corollary \ref{C3.10} shows that the spaces in part (iv) with $s<
1/p$ and also $B^{1/p}_{p,\infty} (\rn)$ do not have the strong truncation property. The case $B^{1/p}_p (\rn)$, $0<p<\infty$, is
again a matter of the Fubini property and Theorem \ref{T3.12}(iii). But this reduction to the one--dimensional case can be extended to
all spaces $B^{1/p}_{p,q} (\rn)$ with $0<p<\infty$ and $0<q \le \infty$. For this purpose one can use the characterization of 
$B^{1/p}_{p,q} (\rn)$ according to \cite[Theorem 2.6.1, (4), p.\,140]{T92} in terms of differences $\Delta^M_{t,k} f$ with $k=1, \ldots, 
n$, where only differences with respect to the $x_k$--directions are involved (a substitute of the Fubini property). Then one argue as
in the one--dimensional case as mentioned in Remark \ref{R3.13} based on the references given there. We do not go into the details.
\end{remark}

For the $F$--spaces one can use the Fubini property according to Proposition \ref{P2.5}(ii) to shift some assertions from $\Fs (\real)
$ to $\Fs (\rn)$, $2 \le n \in \nat$. On the other hand the underlying Theorem \ref{T3.14} for the spaces
$\Fs (\real)$ is less complete than its
counterpart Theorem \ref{T3.12} for the spaces $\Bs (\real)$. We collect the outcome.

\begin{theorem}   \label{T3.18}
Let $n\in \nat$.
\cm
{\em (i)} The spaces $\Fs (\rn)$ with
\begin{\eq}   \label{3.35}
0<p<\infty, \quad 0<q \le \infty \quad \text{and} \quad \sigma^n_{p,q} <s< 1 + \frac{1}{p},
\end{\eq}
$(s \not= \frac{1}{p}$ for $0<p \le 1)$ have the truncation property.
\cm
{\em (ii)} The spaces $\Fs (\rn)$ with
\begin{\eq}   \label{3.36}
0<p<\infty, \quad 0<q \le \infty \quad \text{and} \quad n \max\big( \frac{1}{p}, \frac{1}{q} \big) <s<1
\end{\eq}
have the perfect truncation property.
\cm
{\em (iii)} The spaces $\Fs (\rn)$ with
\begin{\eq}   \label{3.37}
1<p<\infty, \quad 1<q \le \infty \quad \text{and} \quad \max \big( \frac{1}{p}, \frac{1}{q} \big) <s<1
\end{\eq}
have the strong truncation property.
\cm
{\em (iv)} The spaces $\Fs (\rn)$ with
\begin{\eq}   \label{3.38}
0<p<\infty, \quad 0<q \le \infty \quad \text{and} \quad \sigma^n_p <s< \frac{1}{p}
\end{\eq}
do not have the strong truncation property.
\end{theorem}

\begin{remark}   \label{R3.19}
Part (i) is covered by Theorem \ref{T3.3}(ii). Part (ii) follows from Theorem \ref{T2.17}(ii) and \eqref{3.29} in the same way as for
the $B$--spaces in Remark \ref{R3.17} where we added some related comments. Perfect truncation cannot be transferred by the Fubini
property from $\Fs (\real)$ to $\Fs (\rn)$. But  the Fubini property according to Proposition \ref{P2.5}(i) and Theorem 
\ref{T3.14}(ii) reduced to strong truncation prove the above part (iii). Part (iv) is a special case of Corollary \ref{C3.10}.
\end{remark}

\begin{remark}    \label{R3.20}
It has already been observed in \cite{MiVS19} and \cite{VScha20} that the spaces $B^s_p (\rn) =  B^s_{p,p} (\rn)$ with $0 < \frac{1}{p}
<s<1$ have the strong truncation property, whereas the spaces $B^s_p (\rn)$ with $0<s \le \frac{1}{p}$ do not have the strong 
truncation property. The arguments are based on different types of oscillations than above and (at least partly) on the Fubini
property for these spaces. Furthermore, it has shown in \cite{MuN19} that the spaces
\begin{\eq}  \label{3.39}
H^s (\rn) = H^s_2 (\rn) = B^s_2 (\rn) \quad \text{with} \quad \frac{1}{2} <s< \frac{3}{2},
\end{\eq}
have the strong truncation property. This is based on Hilbert space arguments which cannot be extended to the Sobolev spaces
\begin{\eq}   \label{3.40}
H^s_p (\rn) = F^s_{p,2} (\rn), \qquad 1<p<\infty,
\end{\eq}
if $p \not= 2$. On the one hand, \eqref{3.39} is a special case of Theorem \ref{T3.16}(iii), restricted to $B^s_p (\rn)$ with
\begin{\eq}   \label{3.41}
1<p<\infty, \qquad \frac{1}{p} <s< 1 + \frac{1}{p}.
\end{\eq}
On the other hand  it follows from Theorem \ref{T3.18}(iii) that the spaces $H^s_p (\rn) = F^s_{p,2} (\rn)$ with
\begin{\eq}  \label{3.42}
1<p<\infty, \qquad \max \big( \frac{1}{p}, \frac{1}{2} \big) <s<1
\end{\eq}
have the strong truncation property. But there is a gap between \eqref{3.39}, \eqref{3.41} compared with with \eqref{3.42}. Our 
arguments resulting in Theorem \ref{T3.16}(iii) rely on the Fubini property and Theorem \ref{T3.12}(ii), which in turn is based on
Proposition \ref{P2.11}(i) dealing with Faber expansions for $B$--spaces and the rather specific discussion afterwards. One may ask to
employ the $F$--part of Proposition \ref{P2.11}(ii). This might be possible, but has not yet been done and requires surely some
additional efforts. But one take the discussion in Remark \ref{R2.12} as a guide. In particular, \eqref{2.49} suggests to ask whether
the  spaces $H^s_p (\rn) = F^s_{p,2} (\rn)$ have the strong truncation property if
\begin{\eq}  \label{3.43}
\begin{cases}
2 \le p <\infty, &\frac{1}{2} <s< 1 + \frac{1}{p}, \\
1<p<2, & \frac{1}{p} <s< \frac{3}{2}.
\end{cases}
\end{\eq}
The different restrictions \eqref{3.41} for $B^s_p (\rn)$ and \eqref{3.43} for $H^s_p (\rn)$ are quite natural at least in the context
of (one--dimensional) Faber bases and may be found in \cite[Section 3.1.2, pp.\,126--129]{T10}. They, and also the more involved 
counterparts for the spaces $\Fs (\real)$, came out in \cite{T10} as a result of the technicalities used there in connection with 
Haar bases. But somewhat surprisingly it has been discovered later on that they are natural. One may consult \cite[Section 3.5,
pp.\,98--103]{T20} where one finds discussions and corresponding references. In any case it might be a challenging task to find out
whether this remarkable splitting between $B$--spaces on the one hand and $H$-spaces or (more general) $F$--spaces on the 
other hand has a natural counterpart in the connection with perfect truncation (in one dimension) and strong truncation (in higher
dimensions).
\end{remark}

\section{Compositions}    \label{S4}
\subsection{Definitions}    \label{S4.1}
Let
\begin{\eq}    \label{4.1}
 g(t) \quad \text{for $t \in \real$ with $g(0) =0$}
\end{\eq}
be a real and continuous function on $\real$
and let $\As (\rn) \subset L^{\loc}_1 (\rn)$ (locally Lebesgue integrable functions on \rn). Then it makes sense to ask whether
\begin{\eq}   \label{4.2}
g \circ f = g \big( f(\cdot) \big) \in \As (\rn) \quad \text{for any $f\in \As (\rn)$ real}.
\end{\eq}
This so--called {\em composition} (of $f$ with $g$) attracted a lot of attention since decades. We refer the reader to
the survey \cite{BoS11}, the recent publication \cite{BMS20}, the remarkable papers \cite{BMS10}, \cite{BMS14}
and the literature within, including \cite{RuS96} (reflecting the early history of this topic). The most distinguished case of composition is truncation based on $g(t) = |t|$, $t \in \real$. We fixed in 
Definition \ref{D3.1} what is meant by truncation, strong truncation and perfect truncation. One may ask for
an extension of this definition from truncation to composition. The counterpart of the {\em sublinearity} 
\eqref{3.2}, this means
\begin{\eq}    \label{4.3}
\| g \circ f \, | \As (\rn) \| \le c_g \, \| f \, | \As (\rn) \|, \qquad f\in \As (\rn) \ \text{real},
\end{\eq}
is no longer part of the definition \eqref{4.2}. But it is a significant topic to clarify under which 
specifications for $g$ and, in particular, for $\As (\rn)$, the composition \eqref{4.2} has in addition the sublinearity
property \eqref{4.3}. One may consult the above references. Afterwards it  makes sense to ask for 
counterparts of strong truncation and perfect truncation as introduced in Definition \ref{D3.1}, where $L^{\loc}_1 (\rn)$ and
$C(\rn)$ have the same meaning as in the Sections \ref{S2.2} and \ref{S3.1}.

\begin{definition}   \label{D4.1}
Let $g$ be a real and continuous function on $\real$ with $g(0) =0$. Let $n\in \nat$.
\cm
{\em (i)} Let $\As (\rn) \subset L^{\loc}_1 (\rn)$. Then $\As (\rn)$ is said to have the strong $g$--
composition property if $g \circ f \in \As (\rn)$ for any real $f \in \As (\rn)$ and
\begin{\eq}   \label{4.4}
\| g \circ f \, | \As (\rn) \| \sim \| f \, | \As (\rn) \|, \qquad f\in \As (\rn) \ \text{real}.
\end{\eq}
$($equivalences$)$.
\cm
{\em (ii)} Let $\As (\rn) \hra C(\rn)$ $($continuous embedding$)$. Then $\As (\rn)$ is said to have the perfect $g$--composition 
property when any real $f \in C(\rn)$ belongs to $\As (\rn)$ if, and only if, $g \circ f$ belongs to $\As (\rn)$ and
\begin{\eq}   \label{4.5}
\| g \circ f \, | \As (\rn) \| \sim \| f \, | \As (\rn) \|.
\end{\eq}
$($equivalences$)$.
\end{definition}

\begin{remark}   \label{R5.46}
This extends the parts (ii) and (iii) of Definition \ref{D3.1} from $g(t) =|t|$ to the above functions
$g(t)$. We discussed in Remark \ref{R3.2} why the assumption $\As (\rn) \hra C(\rn)$ in connection with
perfect truncation is reasonable. This applies also to the above situations.
\end{remark}

\subsection{Perfect compositions}    \label{S4.2}
A real function $g(t)$ on $\real$ with $g(0) =0$ is said to be a {\em Lipschitz scaling function}  if there 
are two numbers $L_1$, $L_2$ with $0<L_1 \le L_2 <\infty$ such that
\begin{\eq}    \label{4.6}
L_1 (t_2 - t_1) \le g(t_2) - g(t_1) \le L_2 (t_2 - t_1), \quad   - \infty < t_1 \le t_2 <\infty.
\end{\eq}
Then $g(t)$  is a Lipschitz distortion of $t$
and $|g(t)|$ is a Lipschitz distortion of $|t|$. In particular $L_1 |t| \le |g(t)| \le L_2
|t|$, $t\in \real$. Let $0<p \le \infty$. Then
\begin{\eq}   \label{4.7}
L_1 \| f \, | L_p (\rn) \| \le \| g \circ f \, | L_p (\rn) \| \le L_2 \| f \, | L_p (\rn) \|, \quad f\in L_p
(\rn) \ \text{real}
\end{\eq}
and 
\begin{\eq}   \label{4.8}
\| g \circ f \, | L_p (\rn) \| = \| \, |g|\circ f \, | L_p (\rn) \|, \quad f\in L_p (\rn) \ \text{real},
\end{\eq}
follow from the approximation of $f$ by step functions (with respect to cubes) and related limits.

\begin{theorem}   \label{T4.3}
 Let $n\in \nat$. Let $g$ be a Lipschitz scaling function. Then the spaces $\As (\rn)$ with $0<p,q \le \infty$ $(p<\infty$ for 
$F$--spaces$)$,
\begin{\eq}   \label{4.9}
\frac{n}{p} <s<1 \ \text{for $B$--spaces}, \quad n \cdot \max \Big( \frac{1}{p}, \frac{1}{q} \Big) <s<1 \ \text{for $F$--spaces},
\end{\eq}
have both the perfect $g$--composition property and the perfect $|g|$--composition property. 
\end{theorem}

\begin{remark}   \label{R4.4}
This is the extension of the Theorems \ref{T3.16}(ii) and \ref{T3.18}(ii) from perfect truncation to perfect $g$--composition and $|g|$--composition. Instead of \eqref{3.29} and its $n$--dimensional generalization one relies on
\begin{\eq}   \label{4.10}
[g \circ f ]_{j,m} \sim [ |g|\circ f ]_{j,m} \sim [f]_{j,m}, \qquad j\in \no, \quad m \in \zn.
\end{\eq}
and Theorem \ref{T2.17} using in addition \eqref{4.6}, \eqref{4.7}.
\end{remark}

The perfect truncation property and the strong truncation property according to Theorem \ref{T3.12}(ii) rely on rather specific assertions
about Faber bases on \real. It is not clear how to extend  the related arguments from $|t|$ to $g(t)$ and $|g|(t) = |g(t)|$ underlying
Definition \ref{D4.1}. But as a substitute one can try to shift assertions for strong $g$--compositions and perfect $g$--compositions
with the help of Lipschitz diffeomorphisms for some spaces $\Bs (\real)$ to their truncation counterparts. This works quite well but
requires some efforts. First we describe the underlying ingredients which may be of self--contained interest.

\begin{proposition}   \label{P4.5}
Let           
\begin{\eq}   \label{4.11}
1<p \le \infty, \quad 0 < q \le \infty, \quad \text{and} \quad  0<s< 1 + \frac{1}{p}.
\end{\eq}
Let $g\in B^{1+ \frac{1}{p}}_{p,1} (\real)$ be real and $g(0) =0$. Then
\begin{\eq}   \label{4.12}
\| g \circ f \, | \Bs (\rn) \| \le c \, \|g \, | B^{1+ \frac{1}{p}}_{p,1} (\real) \| \cdot \|f \, | \Bs (\rn) \|, \quad f\in \Bs (\rn)
\ \text{real}.
\end{\eq}
\end{proposition}

\begin{remark}   \label{R4.6}
Recall that
\begin{\eq}   \label{4.13}
B^{1+ \frac{1}{p}}_{p,1} (\real) \hra B^1_{\infty,1} (\real) \hra C^1 (\real), \qquad 0<p \le \infty,
\end{\eq}
where $C^1 (\real)$ collects all functions $u(t)$, $t\in \real$, such that both $u$ and its classical  first derivative $u'$ are
bounded and continuous on \real. This shows that $g$ and $g'$ fit in the above scheme. The proof in  \cite[Proposition 5.48]{Tri21}
relies on some non--linear interpolation and the following crucial observation. Let
\begin{\eq}   \label{4.14}
1<p<\infty, \quad 1 \le q \le \infty \quad \text{and} \quad 1<s< 1 + \frac{1}{p}.
\end{\eq}
Then
\begin{\eq}   \label{4.15}
\| g \circ f \, | \Bs (\rn) \| \le c \Big( \| g' \, | L_\infty (\real) \| + \| g' \, | \dot{B}^{1/p}_{p, \infty} (\real) \| \Big) \,
\|f \, | \Bs (\rn) \|
\end{\eq}
is the main assertion in \cite{Kat00}, $f\in \Bs (\rn)$ real, in the reformulation according to \cite[Proposition 3.1]{BMS20}. Here
$u \in \dot{B}^\sigma_{p,\infty} (\real)$ for $0<p<\infty$ and $0<\sigma <1$ means that
\begin{\eq}   \label{4.16}
\| u\, | \dot{B}^\sigma_{p,\infty} (\real) \| = \sup_{h>0} h^{-\sigma} \| u(\cdot + h) - u(\cdot)\, | L_p (\real)\|
\end{\eq}
is finite. An extension of \eqref{4.14} to $s = 1 + \frac{1}{p}$ in the context of \eqref{4.15} goes back to \cite{Mou18}. There one
finds in term of the above reformulation
\begin{\eq}    \label{4.17}
\| g \circ f \, | B^{1+ \frac{1}{p}}_{p,\infty} (\rn) \| \le c\, \| g \, | B^{1+ \frac{1}{p}}_{p,1} (\real) \| \cdot \|f \, | 
B^{1 + \frac{1}{p}}_{p,1} (\rn) \|, \quad f\in B^{1+ \frac{1}{p}}_{p,1} (\rn) \ \text{real}.
\end{\eq}
This assertion and non--linear interpolation are used in \cite{Mou18} to recover \eqref{4.15} with $p,q$ as in \eqref{4.14} and now
$0<s<1 + \frac{1}{p}$ as in \eqref{4.11}.
\end{remark}

Let $g$ be a Lipschitz scaling function  according to  \eqref{4.6} and let $g^{-1}$ be its inverse, $g(t) = \tau$ if, and only if, 
$g^{-1} (\tau) =t$ where $t\in \real$ and $\tau \in \real$. Let $1<p< \infty$ and let, in addition, $g' \in \dot{B}^{1/p}_{p,\infty} 
(\real)$ in the understanding of \eqref{4.16}. Then it follows from some elementary calculations as detailed in \cite[Proposition 
5.58]{Tri21} that $g^{-1}$ is again a Lipschitz scaling function and $\big( g^{-1} \big)' \in \dot{B}^{1/p}_{p,\infty} (\real)$. Now
one apply \eqref{4.15} both to $g$ and $g^{-1}$ with the following outcome.

\begin{proposition}   \label{P4.7}
Let $n \in \nat$,
\begin{\eq}   \label{4.18}
1<p \le \infty, \qquad 0<q \le \infty \quad \text{and} \quad 0<s<1 + \frac{1}{p}.
\end{\eq}
Let $g$ be a Lipschitz scaling function according to \eqref{4.6}. Let, in addition $g' \in \dot{B}^{1/p}_{p, \infty} (\real)$ if $s \ge
1$. Then  any real function $f$ belongs to
$\Bs (\rn)$ if, and only if, $ g \circ f \in \Bs (\rn)$. Furthermore,
\begin{\eq}   \label{4.19}
\| g \circ f \, | \Bs (\rn) \| \sim \| f \, | \Bs (\rn) \|, \qquad f\in \Bs (\rn) \quad \text{real}.
\end{\eq}
\end{proposition}

\begin{remark}   \label{R4.8}
This follows from Proposition \ref{P4.5} in the version of \eqref{4.15}, the above comments and $g^{-1} \circ g = \id$. If $s<1$ then
the spaces $\Bs (\rn)$ can be characterized by first differences similarly as in \eqref{4.16}. Now \eqref{4.19} follows by elementary
arguments.
\end{remark}

\begin{remark}   \label{R4.9}
The above assertion fits in the scheme of perfect $g$--composition according to Definition \ref{D4.1}(ii). But it does not say anything
about the more interesting strong or perfect $|g|$--composition. One may \eqref{4.19} better call a non-linear Lipschitz distortion
of the real part of the space $\Bs (\rn)$ with \eqref{4.18}.
\end{remark}

Let the strong $g$--composition and the perfect $g$--composition be as introduced in Definition \ref{D4.1}. We wish to complement
Theorem \ref{T3.12}(ii) replacing  $|t|$ by $|g(t)|$ where $g$ is a Lipschitz scaling function according to \eqref{4.6}. For this
purpose we combine this theorem with Proposition \ref{P4.7}. 

\begin{theorem}    \label{T4.10}
Let $g$ be a Lipschitz scaling function according to \eqref{4.6}. Let
\begin{\eq}   \label{4.20}
1<p \le \infty, \quad 0<q \le \infty \quad \text{and} \quad \frac{1}{p} <s< 1+ \frac{1}{p}.
\end{\eq}
Let, in addition, $g' \in \dot{B}^{1/p}_{p,\infty} (\real)$ if $s \ge 1$ in the understanding of \eqref{4.16}. 
Then the spaces $\Bs (\real)$ have both the perfect $g$--composition property and the perfect $|g|$--composition property.
\end{theorem}

\begin{remark}   \label{R4.11}
The perfect $g$--composition property is already covered by Proposition \ref{P4.7}. It follows from the properties of the Lipschitz
scaling function $g$ that
\begin{\eq} \label{4.21}
(|g|\circ f)(x) = \big| g \big( f(x) \big) \big|, \qquad f \in \Bs (\real) \quad \text{real}, \quad x\in \real.
\end{\eq}
Otherwise one can combine Proposition \ref{P4.7} and Theorem \ref{T3.12}(ii). But one can argue more directly, following the discussion
in \eqref{2.60} and afterwards, combined with the observation that $g \circ f$ respects the sign of $f$.
\end{remark}

\subsection{Strong compositions}    \label{S4.3}
Similarly as for truncations we reduce strong compositions for spaces $\As (\rn)$ in \rn, $n \ge 2$, to perfect compositions for the
related spaces $\As (\real)$ on $\real$ and the Fubini property according to Proposition \ref{P2.5}. Compositions have the same
meaning as in \eqref{4.2} and Definition \ref{D4.1}. Lipschitz scaling functions have been introduced at the beginning of Section
\ref{S4.2}. Let again $B^s_p (\rn) = B^s_{p,p} (\rn)$, $s\in \real$ and $0<p \le \infty$.

\begin{theorem}    \label{T4.12}
Let $n\in \nat$. 
\cm
{\em (i)} Let $g$ be a Lipschitz scaling function. Then the spaces $\Fs (\rn)$ with
\begin{\eq}   \label{4.22}
1<p<\infty, \quad 1<q \le \infty \quad \text{and} \quad \max \big( \frac{1}{p}, \frac{1}{q} \big) <s<1
\end{\eq}
have both the strong $g$--composition property and the strong $|g|$--composition property.
\cm
{\em (ii)} Let $g$ be a Lipschitz scaling function. Let
\begin{\eq}   \label{4.23}
1<p \le \infty \quad \text{and} \quad \frac{1}{p} < s < 1 + \frac{1}{p}.
\end{\eq}
Let, in addition, $g' \in \dot{B}^{1/p}_{p,\infty} (\real)$ if $s \ge 1$. Then the spaces $B^s_p (\rn)$ have both the strong 
$g$--composition property and the strong $|g|$--composition property.
\end{theorem}

\begin{remark}   \label{R4.13}
Part (i) follows from Theorem \ref{T4.3} (with $n=1$) and the Fubini property according to Proposition \ref{P2.5}. Part (ii) is a
consequence of Theorem \ref{T4.10} and again the Fubini property.
\end{remark}

\begin{remark}   \label{R4.14}
It follows from the above theorem that the Sobolev spaces
\begin{\eq}  \label{4.24}
H^s_p (\rn) = F^s_{p,2} (\rn), \qquad 1<p<\infty,
\end{\eq}
have the strong $g$--composition property and the strong $|g|$--composition property if
\begin{\eq}   \label{4.25}
1<p<\infty \quad \text{and} \quad \max \big( \frac{1}{p}, \frac{1}{2} \big) <s <1.
\end{\eq}
This is the same as for the strong truncation property as we discussed  in Remark \ref{R3.20}. One can ask as there whether 
\eqref{3.43} is a more natural restriction  for the strong composition property for $H^s_p (\rn)$, again in contrast to the natural 
restriction \eqref{4.23} for the spaces $B^s_p (\rn)$.
\end{remark}

%\newpage


\begin{thebibliography}{i}

\bibitem{BMS10} G. Bourdaud, M. Moussai, W. Sickel. Composition operators on Lizorkin--Triebel spaces.
J. Funct. Anal. {\bfseries 259} (2010), 1098--1128.

\bibitem{BMS14} G. Bourdaud, M. Moussai, W. Sickel. Composition operators acting on Besov spaces on the
real line. Ann. Mat. Pura Appl. {\bfseries 193} (2014), 1519--1554.

\bibitem{BMS20} G. Bourdaud, M. Moussai, W. Sickel. A necessary condition for composition in Besov spaces. Compl. Var. Elliptic Equ. {\bfseries 65} (2020), 22--39. 

\bibitem{BoS11} G. Bourdaud, W. Sickel. Composition operators on function spaces with fractional order 
of smoothness. RIMS K\^{o}ky\^{u}roku Bessatsu {\bfseries B26}, Res. Inst. Math. Sci. (RIMS), Kyoto, 2011,
93--132.

\bibitem{GiT77} D. Gilbarg, N.S. Trudinger. Elliptic partial differential equations of second order. Springer, Berlin, 1977.

\bibitem{Kat00} D. Kateb. Fonctions qui operent sur les espaces de Besov. Proc. Amer. Math. Soc. 
{\bfseries 128} (2000), 735--743.

\bibitem{MiVS19} P. Mironescu, J. Van Schaftingen. Lifting in compact covering spaces for fractional Sobolev mappings. arXiv:
1907.01373 (2019).

\bibitem{Mou18} M. Moussai. Composition operators on Besov spaces in the limiting case $s= 1 + 
\frac{1}{p}$. Studia Math. {\bfseries 241} (2018), 1--15. 

\bibitem{MuN19} R. Musina, A.I. Nazarov. A note on truncations in fractional Sobolev spaces. Bull. Math. Sci. {\bfseries 9}
(2019), 1950001 (7 pages). 


\bibitem{RuS96} T. Runst, W. Sickel. Sobolev spaces of fractional order, Nemytskij operators, and nonlinear partial differential
equations. W. de Gruyter, Berlin, 1996.

\bibitem{Tri89} H. Triebel. Local approximation spaces. Z. Anal. Anwend. {\bfseries 8} (1989), 261--288.

\bibitem{T92} H. Triebel. Theory of function spaces II. Birkh\"{a}user, Basel, 1992. 

\bibitem{T01} H. Triebel. The structure of functions. Birkh\"{a}user, Basel, 2001.

\bibitem{T10} H. Triebel. Bases in function spaces, sampling, discrepancy, numerical integration. European Math. Soc. Publishing House, Z\"{u}rich, 2010.

\bibitem{T12} H. Triebel. Faber systems and their use in sampling, discrepancy, numerical integration. European Math. Soc. Publishing House, Z\"{u}rich, 2012.

\bibitem{T15} H. Triebel. Tempered homogeneous function spaces. European Math. Soc. Publishing House, Z\"{u}rich, 2015.

\bibitem{T20} H. Triebel. Theory of function spaces IV. Birkh\"{a}user, Basel, 2020. 

\bibitem{Tri21} H. Triebel. Quarkonial analysis aand truncations. Informal Notes, Jena, 2020.

\bibitem{VScha20} J. Van Schaftingen. Reverse superposition estimates in Sobolev spaces. arXiv: 
2007.01742 (2020).

\bibitem{Zie89} W.P. Ziemer. Weakly differentiable functions. Springer, New York, 1989.



\end{thebibliography}
\end{document}